\documentclass[a4paper,12pt]{article}

\usepackage[latin1]{inputenc}
\usepackage{natbib}
\usepackage{amsmath}
\usepackage{amsthm}
\usepackage{amssymb} 
\newcommand{\rea}{\mathbb{R}}

\usepackage{stmaryrd} 
\usepackage{ctable}
\usepackage{MnSymbol}
\usepackage{natbib}
\usepackage{mathtools}
\usepackage{endnotes} 
\usepackage{setspace}

\newtheorem{Def}{Definition}
\newtheorem{Theo}{Theorem}
\newtheorem{cor}{Corollary}
\newtheorem{rem}{Remark}

\begin{document}

\title {Bipolar Fuzzy Integrals}

\author{Salvatore Greco \\ \small{\textit{Department of Economics and Business, University of Catania, Corso Italia 55, 95129 Catania, Italy}}\\\small{E-mail: {salgreco@unict.it}}
\\
Fabio Rindone\footnote{Corresponding author: telephone +390957537733, fax +390957537957} \\ \small{\textit{Department of Economics and Business, University of Catania, Corso Italia 55, 95128 Catania, Italy}}\\\small{E-mail: {frindone@unict.it}}
}

\date{}
\maketitle

\begin{abstract}
In decision analysis and especially in multiple criteria decision analysis, several non additive integrals have been introduced in the last sixty years. 
Among them, we remember the Choquet integral, the Shilkret integral and the Sugeno integral. 
Recently, the bipolar Choquet integral has been proposed for the case in which the underlying scale is bipolar. 
In this paper we propose the bipolar Shilkret integral and the bipolar Sugeno integral. 
Moreover, we provide an axiomatic characterization of all these three bipolar fuzzy integrals.  
\end{abstract}


\section{Introduction}
In decision analysis and especially in multiple criteria decision analysis, several non additive integrals have been introduced in the last sixty years \cite{figueira2005multiple,grabisch1996application, grabisch2005fuzzy}.  
Among them, we remember the Choquet integral \cite{choquet1953theory}, the Shilkret integral \cite{shilkret1971maxitive} and the Sugeno integral \cite{sugeno1974theory}.        
Recently the bipolar Choquet integral \cite{grabisch2005biI,grabisch2005bi,greco2002bipolar} has been proposed for the case in which the underlying scale is bipolar. 
A further generalization is that of level dependent integrals, which has lead to the definition of the level dependent Choquet integral \cite{greco2011choquet}, the level dependent Shilkret integral \cite{bodjanova2009sugeno}, the level dependent Sugeno integral \cite{mesiar2009level} and the bipolar level dependent Choquet integral \cite{greco2011choquet}.
Very recently, on the basis of a minimal set of axioms, one concept of universal integral giving a common framework to many of the above integrals have been proposed \cite{klement2010universal}. 
In this paper we aim to provide a general framework for the case of bipolar fuzzy integrals, i.e. those integrals whose underlying scale is bipolar. 
For this purpose we propose the definition of bipolar Shilkret integral and bipolar Sugeno integral. 
Then, in order to provide a mathematical characterization of the three mentioned bipolar integrals, we give necessary and sufficient conditions for an aggregation function to be the bipolar Choquet integral or the bipolar Shilkret integral or the bipolar Sugeno integral. 
As we said, the bipolar fuzzy integrals admit a further generalization if 
the fuzzy measure (capacity) with respect to which the integrals are calculated can change from a level to another \cite{greco2011choquet,greco2011linz}. 
For the sake of clarity, we shall remind the characterization of the bipolar Shilkret and Sugeno integral with respect to a level dependent capacity in a forthcoming paper (we wish to remember as such results have just been presented in \cite{greco2011linz}). 
The paper is organized as follows. 
In section \ref{Preliminaries} we give the preliminaries and list some  properties of an aggregation function useful to the characterization of the bipolar fuzzy integrals we shall propose in this paper. 
In section \ref{fuzzyintegrals} we review the definitions and characterizations of the classical Choquet integral, Shilkret integral and Sugeno integral.
In section \ref{bipfuzzy} we give our main results: first we propose the bipolar version of the Shilkret integral and of the Sugeno integral; next we characterize the bipolar Choquet, Shilkret and Sugeno integrals. 
Section \ref{concluding} contains conclusions. 
All the proofs are presented in the appendix. 

\section{Preliminaries}\label{Preliminaries}
Let us consider a set of criteria 
$N=\left\{1,\ldots,n\right\}$ and let 
$(\alpha,\beta)$ 
be any possible interval of $\rea$, i.e any of 
$[-\alpha,\beta], \ ]-\alpha,\beta],\ 
[-\alpha,\beta[,\ 
]-\alpha,\beta[,\ 
]-\infty  ,\beta[, \ 
]-\infty  ,\beta], \ 
]-\alpha, +\infty[,\ 
[-\alpha, +\infty[,\  
]-\infty,+\infty[
.$
An \textit{alternative} can be identified with a score vector $\textbf{\textit{x}}=\left(x_1, \ldots, x_n\right)\in \left(\alpha,\beta\right)^n$, being $x_i$ the evaluation of such an alternative $\textbf{\textit{x}}$ with respect to the $i^{th}$ criterion. 
An alternative $\textbf{\textit{x}}$ dominates another $\textbf{\textit{y}}$ if on each criterion the evaluation of $\textbf{\textit{x}}$ is not smaller than the evaluation of $\textbf{\textit{y}}$, i.e. for all $i\in N$ $x_i\ge y_i$ and in this case we simply write $\textbf{\textit{x}}\ge\textbf{\textit{y}}$.
The indicator function of any $A\subseteq N$ is the function which attains 1 on $A$ and and 0 on $N\setminus A$ and can be identified with the vector $\textbf{1}_A$ whose $i^{th}$ component is equal to 1 if $i \in A$ and 0 otherwise.
\\ 
In general, an aggregation function is a function $G: (\alpha,\beta)^n \rightarrow(\alpha,\beta)$ such that
\begin{enumerate}
\item 
	$G(\alpha, \ldots, \alpha)=\alpha$ if $\alpha \in (\alpha,\beta)$ 
	and 
	$\lim_{x\rightarrow \alpha^+}G(x,\ldots,x)=\alpha$ if $\alpha \notin (\alpha,\beta);$ 
\item 
$G(\beta, \ldots, \beta)=\beta$ if $\beta \in (\alpha,\beta)$ 
and 
$\lim_{x\rightarrow \beta^-}G(x,\ldots,x)=\beta$ if $\beta \notin (\alpha,\beta)$;
\item for all $\textbf{\textit{x}},\textbf{\textit{y}} \in (\alpha,\beta)^n$ such that  
$\textbf{\textit{x}} \ge \textbf{\textit{y}}$, $G(\textbf{\textit{x}}) \ge G(\textbf{\textit{y}})$.
\end{enumerate}

\noindent In this paper we often denote the maximum and the minimum of a set $X$ respectively with $\bigvee X$ and $\bigwedge X$. 
For any two alternatives 
$\textbf{\textit{x}}, \textbf{\textit{y}} \in \left(\alpha,\beta\right)^n$, 
the following definitions hold  
\begin{itemize}
  \item $\textbf{\textit{x}} \wedge \textbf{\textit{y}}$ is the vector whose $i^{th}$ component is  
$\left(x\wedge y\right)_i=\bigwedge\{x_i, y_i\}$  for all $i=1, \ldots, n$ (in case $\textbf{\textit{y}}=\left(h,\ldots,h\right)$ is a constant, then we can write $\textbf{\textit{x}} \wedge h$);
  \item $\textbf{\textit{x}} \vee \textbf{\textit{y}}$ is the vector whose $i^{th}$ component is 
$\left(x\vee y\right)_i=\bigvee\{x_i, y_i\}$  for all $i=1, \ldots, n$ (in case $\textbf{\textit{y}}=\left(h,\ldots,h\right)$ is a constant, then we can write $\textbf{\textit{x}} \vee h$);
  \item $\textbf{\textit{x}}$ and $\textbf{\textit{y}}$ are comonotone (or comonotonic) if $(x_i-x_j)(y_i-y_j) \ge 0$ for all $i,j \in N$;
  \item $\textbf{\textit{x}}$ and $\textbf{\textit{y}}$ are bipolar comonotone if 
$(|x_i|-|x_j|)(|y_i|-|y_j|) \ge 0$ and $x_iy_i\ge 0$, 
for all $i,j \in N$.
\end{itemize}
The following properties of an aggregation function $G:\left(\alpha,\beta\right)^n \rightarrow\left(\alpha,\beta\right)$ are useful to characterize several integrals: 
\begin{itemize}
	\item idempotency: for all $\textbf{a}\in\left(\alpha,\beta\right)^n$ such that $\textbf{a}=(a,\ldots,a)$, $G(\textbf{a})=a$;
	\item homogeneity: for all $\textbf{\textit{x}} \in\left(\alpha,\beta\right)^n$ and $c>0$ such that $c\cdot\textbf{\textit{x}} \in \left(\alpha,\beta\right)^n$,  $G(c\cdot\textbf{\textit{x}})=c\cdot G(\textbf{\textit{x}})$;
	\item stability w.r.t. the minimum: for all 
	$\textbf{\textit{x}} \in\left(\alpha,\beta\right)^n$ and $\gamma \in \left(\alpha,\beta\right)$, $G(\textbf{\textit{x}} \wedge \gamma)=\bigwedge\{G(\textbf{\textit{x}}),\gamma\}$; 
	\item additivity: for all $\textbf{\textit{x}}, \textbf{\textit{y}}\in\left(\alpha,\beta\right)^n$ such that $\textbf{\textit{x}}+\textbf{\textit{y}}\in\left(\alpha,\beta\right)^n$, $G(\textbf{\textit{x}}+\textbf{\textit{y}})=G(\textbf{\textit{x}})+G(\textbf{\textit{y}})$;
	\item maxitivity: for all $\textbf{\textit{x}}, \textbf{\textit{y}}\in\left(\alpha,\beta\right)^n$, with $\alpha\ge 0$,  $G(\textbf{\textit{x}}\vee\textbf{\textit{y}})=\bigvee \{G(\textbf{\textit{x}}),G(\textbf{\textit{y}})\}$; 
	\item minitivity: for all $\textbf{\textit{x}}, \textbf{\textit{y}}\in\left(\alpha,\beta\right)^n$, with $\beta\le 0$,  $G(\textbf{\textit{x}}\wedge\textbf{\textit{y}})=\bigwedge \{G(\textbf{\textit{x}}),G(\textbf{\textit{y}})\}$; 
	\item comonotonic additivity: for all comonotone $\textbf{\textit{x}}, \textbf{\textit{y}}\in\left(\alpha,\beta\right)^n$, $G(\textbf{\textit{x}}+\textbf{\textit{y}})=G(\textbf{\textit{x}})+G(\textbf{\textit{y}})$;
		\item comonotonic maxitivity: for all comonotone $\textbf{\textit{x}}, \textbf{\textit{y}}\in\left(\alpha,\beta\right)^n$,  $G(\textbf{\textit{x}}\vee \textbf{\textit{y}})=\bigvee\{G(\textbf{\textit{x}}),G(\textbf{\textit{y}})\};$
		\item comonotonic minitivity: for all comonotone $\textbf{\textit{x}}, \textbf{\textit{y}}\in\left(\alpha,\beta\right)^n$,  $G(\textbf{\textit{x}}\wedge \textbf{\textit{y}})=\bigwedge\{G(\textbf{\textit{x}}),G(\textbf{\textit{y}})\};$
	\end{itemize}

\section{Fuzzy integrals}\label{fuzzyintegrals}
Let us briefly review the three most famous fuzzy integrals, i.e. the Choquet, Shilkret and Sugeno integrals. 
For each of them we shall discuss the restrictions to be imposed on the scale $\left(\alpha,\beta\right)$. 

\subsection{The Choquet integral}
\begin{Def}
A capacity is function $\mu: 2^N \rightarrow [0,1]$ satisfying the following properties:
\begin{enumerate}
	\item $\mu(\emptyset)=0,\  \mu(N)=1$,
	\item for all $A \subseteq B \subseteq N, \ \mu(A) \le \mu(B)$.
\end{enumerate}
\label{capacity}
\end{Def}

\begin{Def}
The Choquet integral \cite{choquet1953theory} of a vector $\textbf{\textit{x}}=\left(x_1, \ldots, x_n\right)\in \left(\alpha,\beta\right)^n\subseteq \left[0,+\infty\right.\left[\right.^n$ with respect to the capacity $\mu$ is given by
\begin{equation}
Ch(\textbf{\textit{x}},\mu)=\int_{0}^{\infty} \mu\left(\{i\in N: x_i \ge t\}\right) dt.
\label{choquet}
\end{equation}
\end{Def}

\noindent Schmeidler \cite{schmeidler1986integral} extended the above definition to negative values too, moreover he characterized the Choquet integral in terms of comonotonic additivity and idempotency. 

\begin{Def}
The Choquet integral \cite{schmeidler1986integral} of a vector $\textbf{\textit{x}}=\left(x_1, \ldots, x_n\right)\in \left(\alpha, \beta\right)^n$ with respect to the capacity $\mu$ is given by
\begin{equation}
Ch(\textbf{\textit{x}},\mu)=\int_{\min_i x_i}^{\max_i x_i} \mu\left(\{i\in N: x_i \ge t\}\right) dt +\min_i x_i.
\label{choquetschmeidler}
\end{equation}
\end{Def}

\noindent Alternatively \eqref{choquetschmeidler} can be written as 
\begin{equation}
Ch(\textbf{\textit{x}},\mu)=\sum_{i=2}^n\left(x_{(i)}-x_{(i-1)}\right)\cdot\mu\left(\{j\in N: x_j \ge x_{(i)}\}\right)\ +\ x_{(1)}
\label{choquet1}
\end{equation}
being $():N\rightarrow N$ any permutation of indexes such that $x_{(1)}\le\ldots\le x_{(n)}$.

\begin{Theo} \cite{schmeidler1986integral} An aggregation function $G: \left(\alpha,\beta\right)^n\rightarrow\left(\alpha,\beta\right)$ is idempotent and comonotone additive if and only if there exists a  capacity $\mu$ such that, for all $\textbf{\textit{x}}\in \left(\alpha,\beta\right)^n$,
$$G(\textbf{\textit{x}})=Ch(\textbf{\textit{x}},\mu).$$
\end{Theo}

\subsection{The Shilkret integral}
 
\begin{Def}
The  Shilkret integral \cite{shilkret1971maxitive} of a vector $\textbf{\textit{x}}=\left(x_1, \ldots, x_n\right)\in \left(\alpha,\beta\right)^n\subseteq\left[0,+\infty\right.\left[\right.^n$ with respect to the capacity $\mu$ is given by
\begin{equation}
Sh(\textbf{\textit{x}},\mu)=\bigvee_{i \in N}\left\{x_i\cdot\mu(\{j \in N: x_j \ge x_i\}\right\}.
\label{shilkret}
\end{equation}
\end{Def}
\noindent An original result of this paper is the characterization of the Shilkret integral in terms of idempotency, comonotonic maxitivity and homogeneity.
\begin{Theo} Suppose that $\alpha\ge 0$, then an aggregation function $G: \left(\alpha,\beta\right)^n\rightarrow\left(\alpha,\beta\right)$ is idempotent, comonotone maxitive and homogeneous if and only if there exists a capacity $\mu$ on $N$ such that, for all $\textbf{\textit{x}}\in \left(\alpha,\beta\right)^n$,
$$G(\textbf{\textit{x}})=Sh(\textbf{\textit{x}},\mu). $$
\label{theo:shilkret}
\end{Theo}

\begin{rem}\label{rem1}
Given to the importance we give to this result, we shall present a direct proof in the appendix. 
Alternatively, theorem \ref{theo:shilkret} can be elicited as corollary of another theorem we shall present in the next section.
\end{rem} 

\noindent Although in \cite{shilkret1971maxitive} the Shilkret integral was formulated for nonnegative functions, however \eqref{shilkret} works also for a generic $\textbf{\textit{x}}\in\left(\alpha,\beta\right)^n\subseteq \rea^n$. 
But, in our opinion, if we allow for negative values too, the essence of the Shilkret integral
is lost. 
Let us stress this point with some examples. 
Suppose that an alternative is strongly negatively evaluated on each criterion except on the last, where it has a low nonnegative evaluation, e.g. 
$\textbf{\textit{x}}=\left(-100,-100,-100,1\right)$. 
By applying \eqref{shilkret}, 
$Sh\left(\textbf{\textit{x}},\mu\right)=\mu\left(\left\{4\right\}\right)$, for every capacity $\mu$.
Thus, the negative evaluations and the weights that the capacity assigns to the relative criteria with respect to which these negative evaluations are given, are ininfluent on the evaluation of $\textbf{\textit{x}}$. 
In general, if for a given alternative $\textbf{\textit{x}}$ we have simultaneously negative and positive evaluations on the various criteria, the negative ones are ininfluent and the Shilkret integral of $\textbf{\textit{x}}$ coincides with the Shilkret integral of 
$\textbf{\textit{x}}\vee 0$.
In the case of $\textbf{\textit{x}}\in {]-\infty,0[}^n$ it is straightforward noting that 
$Sh\left(\textbf{\textit{x}},\mu\right)=\left(\max_{i\in N}x_i\right)\cdot \mu\left(\left\{j\in N\ |\ x_j\ge\max_{i\in N}x_i \right\}\right)$.
Again, we note how for all capacities only the maximum evaluation of $\textbf{\textit{x}}$ matters. 
For vectors with non-positive evaluation on each criterion, the logic of the Shilkret integral can be recovered if in the \eqref{shilkret} we substitute the maximum with the minimum and $\ge$ with $\le$.

\begin{Def}
The  negative Shilkret integral of a vector $\textbf{\textit{x}}=\left(x_1, \ldots, x_n\right)\in \left(\alpha,\beta\right)^n\subseteq {]-\infty, 0]}^n$ with respect to the capacity $\mu$ is given by
\begin{equation}
Sh^-(\textbf{\textit{x}},\mu)=\bigwedge_{i \in N}\left\{x_i\cdot\mu(\{j \in N: x_j \le x_i\}\right\}=
-\bigvee_{i \in N}\left\{-x_i\cdot\mu(\{j \in N: -x_j \ge -x_i\}\right\}=-Sh(-\textbf{\textit{x}},\mu).
\label{shilkretnegative}
\end{equation}
\end{Def}
\noindent Obviously, from theorem \ref{theo:shilkret}, the characterization of the negative Shilkret integral is in terms of idempotency, comonotonic minitivity and homogeneity.
\begin{cor} Suppose that $\beta\le 0$, then an aggregation function $G: \left(\alpha,\beta\right)^n\rightarrow\left(\alpha,\beta\right)$ is idempotent, comonotone minitive and homogeneous if and only if there exists a capacity $\mu$ on $N$ such that, for all $\textbf{\textit{x}}\in \left(\alpha,\beta\right)^n$,
$$G(\textbf{\textit{x}})=Sh^-(\textbf{\textit{x}},\mu). $$
\label{theo:shilkretnegative}
\end{cor}

\noindent So far, we have a Shilkret integral for alternatives with all non-negative evaluations and one for alternatives with all non-positive evaluations. 
To obtain a suitable definition of the Shilkret integral for the mixed case we propose two different approach.
In the first approach we define a \textit{symmetric Shilkret integral} by applying a logic à la \v{S}ipo\v{s} \cite{sipos1979integral}, i.e. for all $\textbf{\textit{x}}\in\left(\alpha,\beta\right)$
\begin{equation}
\textnormal{\v{Sh}}\left(\textbf{\textit{x}},\mu\right)=Sh(\textbf{\textit{x}}\vee 0,\mu) + Sh^-(\textbf{\textit{x}}\wedge 0,\mu).
\label{symmetricshilkret}
\end{equation}
Note that the \eqref{symmetricshilkret} is called symmetric since 
$\textnormal{\v{Sh}}\left(\textbf{\textit{x}},\mu\right)=-\textnormal{\v{Sh}}\left(-\textbf{\textit{x}},\mu\right)$. 
A second, more general, approach will be to define a \textit{bipolar Shilkret integral} (see next section). 
This would be used directly for the bipolar scale, while restricted on $\rea^+$ and on $\rea^-$ it would coincide respectively with the Shilkret integral and the negative Shilkret integral.

\subsection{The Sugeno integral}

\begin{Def}
A measure on $N$ with a scale $(\alpha,\beta)$ is any function $\nu:2^N\rightarrow (\alpha,\beta)$ such that:
\begin{enumerate}
	\item $\nu(\emptyset)=\alpha,\  \nu(N)=\beta$,
	\item for all $A \subseteq B \subseteq N,\  \nu(A) \le \nu(B)$.
\end{enumerate}
\end{Def}

\begin{Def}
The Sugeno integral \cite{sugeno1974theory} of a vector $\textbf{\textit{x}}=\left(x_1, \ldots, x_n\right)\in\left(\alpha ,\beta\right)^n$ with respect to the measure $\nu$ on $N$ with scale $(\alpha,\beta)$ is given by
\begin{equation}
Su(\textbf{\textit{x}},\nu)=\bigvee_{i\in N}\bigwedge\left\{x_i, \nu\left(\left\{j\in N \ |\ x_j\ge x_i\right\}\right)\right\}.
\label{eq:sugeno}
\end{equation}
\label{sugeno}
\end{Def}
\noindent Alternatively the Sugeno integral can be written as 
\begin{equation}
Su(\textbf{\textit{x}},\nu)=\bigvee_{A \subseteq N}\bigwedge\left\{\nu(A), \bigwedge_{i\in A}x_i\right\}.
\label{eq:sugeno1}
\end{equation}
Next theorem gives necessary and sufficient conditions to be an aggregation function the Sugeno integral.
\begin{Theo} \cite{marichal2001axiomatic} An aggregation function $G: (\alpha ,\beta)^n\rightarrow(\alpha ,\beta)$ is idempotent, comonotone maxitive and stable with respect to the minimum if and only if there exists a  measure $\nu$ on $N$ with a scale $(\alpha,\beta)$ such that, for all $\textbf{\textit{x}}\in (\alpha ,\beta)^n$,
$$G(\textbf{\textit{x}})=Su(\textbf{\textit{x}},\nu). $$
\end{Theo}

\noindent Let us observe that the definition of the Sugeno integral only imposes that the $x_i$ and the $\nu(A)$ are measured on the same (possible only ordinal) scale $\left(\alpha,\beta\right)$. 
Suppose that 
$\mu:2^N\rightarrow [0,1]$, 
is a capacity and  
$\textbf{\textit{x}}\in {[-1,1]}^n$ 
is a vector evaluated on each criterion on the symmetric scale $[-1,1]$, the \textit{symmetric Sugeno integral} \cite{grabisch2003symmetric} of $\textbf{\textit{x}}$ is defined as
\begin{equation}
\textnormal{\v{Su}}\left(\textbf{\textit{x}},\mu\right)=Su(\textbf{\textit{x}}\vee 0,\mu) - Su((-\textbf{\textit{x}})\vee 0,\mu).
\label{symmetricsugeno}
\end{equation}
In \eqref{symmetricsugeno}, as before in \eqref{symmetricshilkret}, symmetric means that  
$\textnormal{\v{Su}}\left(\textbf{\textit{x}},\mu\right)=-\textnormal{\v{Su}}\left(-\textbf{\textit{x}},\mu\right)$. \\
Clearly if $x_i\ge 0$ for all $i\in N$, $\textnormal{\v{Su}}\left(\textbf{\textit{x}},\mu\right)=Su(\textbf{\textit{x}},\mu)$, while if 
$x_i\le 0$ for all $i\in N$,
\begin{equation}
\textnormal{\v{Su}}\left(\textbf{\textit{x}},\mu\right)=\bigwedge_{i\in N}\bigvee\left\{x_i, -\nu\left(\left\{j\in N \ |\ x_j\le x_i\right\}\right)\right\}.
\label{eq:sugenonegative}
\end{equation}
\eqref{eq:sugenonegative} can be considered as a definition of a negative Sugeno integral, for the case in which $\textbf{\textit{x}}$ is negatively evaluated on each criterion.
In the next section we shall propose a more general approach, defining a \textit{bipolar Sugeno integral},  which restricted on $\rea^+$ and on $\rea^-$ coincides respectively with the \eqref{sugeno} and the \eqref{eq:sugenonegative}.

\section{Bipolar fuzzy integrals on the scale [-1,1]}\label{bipfuzzy}
The present work is devoted to the study of bipolar fuzzy integrals, i.e. those integrals useful when the scale underlying the alternatives evaluation is bipolar. 
By the sake of simplicity, trough this section we shall adopt the bipolar scale $[-1,1]$ to present our results. 
However, without loss of the generality, they can be extended to every other symmetric interval of $\rea$, i.e. any of  
$[-\alpha,\alpha], 
]-\alpha,\alpha[,\ 
]-\infty  ,+\infty [$,  
where $\alpha\in\rea^+.$
\\ 
Let us consider the set $\mathcal Q=\left\{\left(A,B\right)\in 2^N\times 2^N\ :\ A\cap B=\emptyset\right\}$ of all disjoint pairs of subsets of $N$. 
With respect to the binary relation $(A,B)\precsim(C,D)$ iff $A\subseteq C$ and $B\supseteq D$, $\mathcal Q$ is a lattice, i.e. a partial ordered set in which any two elements have a unique supremum, 
$(A,B)\vee(C,D)=\left(A\cup C,B\cap D\right),$ and a unique infimum, 
$(A,B)\wedge(C,D)=\left(A\cap C,B\cup D\right)$. 
For all $(A,B), (C,D)\in \mathcal Q$ if $A\subseteq C$ and $B\subseteq D$, we simply write $(A,B)\subseteq (C,D)$.  
For all $(A,B)\in\mathcal Q$ the indicator function $1_{(A,B)}:N\rightarrow \{-1,0,1\}$ is the function which attains 1 on $A$, -1 on $B$ and 0 on $\left(A\cup B\right)^c$. 
Such a function can be identified with the vector $\textbf{1}_{(A,B)}$ whose $i^{th}$ component is equal to 1 if $i \in A$, is equal to $-1$ if $i \in B$ and is equal to 0 otherwise.\\
The \textit{symmetric maximum} of two elements - introduced and discussed in \cite{grabisch2003symmetric,grabisch2004mobius} - is defined by the following binary operation:
\begin{equation}
a\varovee b =\left\{
\begin{array}{ll}
	-\left(|a|\vee |b|\right) & \textnormal{if } b\neq -a \textnormal{ and either } |a|\vee |b|=-a \textnormal{ or } =-b
	\\
	0 & \textnormal{if } b=-a
	\\
	|a|\vee |b| & \textnormal{else.}
\end{array}
\right.
\label{eq:sym-max}
\nonumber
\end{equation}
In \cite{mesiar2010discrete} it has been showed that, on the domain $[-1,1]$, the symmetric maximum coincides with two recent symmetric extensions of the Choquet integral, the \textit{balancing Choquet integral} and the \textit{fusion Choquet integral}, when they are computed with respect to the strongest capacity (i.e. the capacity $\nu:2^N\rightarrow\left[0,1\right]$ which attains zero on the empty set and one elsewhere). 
However, the symmetric maximum of a set $X$ cannot be defined, being $\ovee$ non associative; e.g, suppose that $X=\left\{3,-3,2\right\}$, then $\left(3\varovee -3\right)\varovee 2=2$ or $3\varovee \left(-3\varovee 2\right)=0$, depending on the order.
Several possible extensions of the symmetric maximum for dimension $n, n>2$ have been proposed (see \cite{grabisch2004mobius,grabisch2009aggregation} and also the relative discussion in \cite{mesiar2010discrete}). 
One of these extensions is based on the splitting rule applied to the maximum and to the minimum as described in the following. 
Given $X=\left\{x_1,\ldots,x_m\right\}\subseteq\rea$, 
the \textit{bipolar maximum} of $X$, shortly $\bigvee^b X$, is defined in this manner: 
if there exists an element $x_k\in X$ such that $|x_k|>|x_j|\ \forall j:x_j\neq x_k$ then $\bigvee^b X=x_k$; 
otherwise $\bigvee^b X=0$. 
Clearly, the bipolar maximum is related to the symmetric maximum by means of 
\begin{equation}
{\bigvee}^b X={\bigvee_i^m}\displaystyle^b x_i=\left( {\bigvee_i^m} x_i\right) \ovee \left({\bigwedge_i^m} x_i\right).
\label{eq:sym-max1}
\end{equation}
The following definitions are closely related to the above discussion.
\begin{Def}
Given $X=\left\{x_1,\ldots,x_m\right\}\subseteq\rea$, the \textit{positive bipolar maximum} of $X$, shortly $\bigvee^{b^+} X$, is the element with the greatest absolute value, with the convention that, in the case of two different opposite elements with this property, we choose the non-negative.
\label{def:posbipmax}
\end{Def}
\begin{Def}
Given $X=\left\{x_1,\ldots,x_m\right\}\subseteq\rea$, the \textit{negative bipolar maximum} of $X$, shortly $\bigvee^{b^-} X$, is the element with the greatest absolute value, with the convention that, in the case of two different opposite elements with this property, we choose the non-positive.
\label{def:negbipmax}
\end{Def}
\noindent Following these definitions, if $X=\left\{9, -9, 7, -3\right\}$ thus, 
$\bigvee^{b}X=0$, 
$\bigvee^{b^+}X =9$ and 
$\bigvee^{b^-} X=-9$. 
Clearly the three operators just defined are linked by means of the relation: 
$\bigvee^b X=\bigvee^b\left\{\bigvee^{b^+}{X},\bigvee^{b^-}{X}\right\}$.\\
Given the vectors $\textbf{\textit{x}}^1,\ldots , \textbf{\textit{x}}^k\in\left[-1.1\right]^n$ with $K=\{1,\ldots,k\}$, $\underset{j\in K}{\bigvee^b} \textbf{\textit{x}}_j$ is the vector whose $i^{th}$ component is 
$\bigvee^{b}\{x^1_i,\ldots,x^k_i\}$  for all $i=1, \ldots, n$
and 
$\underset{j\in K}{\bigwedge^b} \textbf{\textit{x}}_j$ is the vector whose $i^{th}$ component is 
$\bigwedge^{b}\{x^1_i,\ldots,x^k_i\}$  for all $i=1, \ldots, n$;\\
The following properties of an aggregation function $G:\left[-1,1\right]^n \rightarrow\left[-1,1\right]$ are useful to characterize several bipolar integrals. 
\begin{itemize}
	\item bipolar comonotonic additivity: for all bipolar comonotone $\textbf{\textit{x}}, \textbf{\textit{y}}\in\left[-1,1\right]^n$,   $$G(\textbf{\textit{x}}+\textbf{\textit{y}})=G(\textbf{\textit{x}})+G(\textbf{\textit{y}});$$
	\item bipolar stability of the sign: for all $r,s \in ]0,1]$ and for all $(A,B)\in \mathcal Q$, 
	$$G(r \textbf{1}_{A,B})G(s \textbf{1}_{A,B})>0 \qquad \textnormal{or} \qquad G(r \textbf{1}_{A,B})=G(s \textbf{1}_{A,B})=0,$$ 
	i.e., in simple words, $G(r \textbf{1}_{(A,B)})$ and $G(s \textbf{1}_{(A,B)})$ have the same sign;
\item bipolar stability with respect to the minimum: for all $r,s \in ]0,1]$ such that $r>s$, and for all $(A,B)\in \mathcal Q$,  
$|G(r \textbf{1}_{(A,B)})|\geq|G(s \textbf{1}_{(A,B)})|$ and, moreover, 
$$\textnormal{if}\qquad |G(r \textbf{1}_{(A,B)})|>|G(s \textbf{1}_{(A,B)})|\qquad \textnormal{then}\qquad |G(s \textbf{1}_{(A,B)})|=s.$$
\end{itemize}

\subsection{A specific property: bipolar comonotone maxitivity}
With a slight abuse of notation we extend the relation of set inclusion to $\mathcal Q$, by defining $(A,B)\subseteq (C,D)$ 
if and only if $A \subseteq C$ and $B\subseteq D$, for all $(A,B), (C,D) \in \mathcal Q$.
Let us suppose to have $k$ different levels $l_1,\ldots,l_k\in\rea$ with
$0<l_1<l_2<\ldots<l_k\le 1$ and a sequence 
$\left\{\left(A_i , B_i\right)\right\}_{i=1,\ldots,k}$ such that
$\left(A_i , B_i\right)\in\mathcal Q$  for all $i=1,\ldots,k$ and 
$\left(A_{i+1} , B_{i+1}\right)\subseteq \left(A_i , B_i\right)$ 
for all $i=1,\ldots,k-1.$ 
The vectors  
$l_i\cdot \textbf{1}_{\left(A_i , B_i\right)}$, $i=1,\ldots,k$ 
are bipolar comonotonic and, moreover, by ordering them with respect to the level $l_i$, then in the vector $l_i\cdot \textbf{1}_{\left(A_i , B_i\right)}$, for each component the elements under the level $l_i$ are the opposite of that under the level $-l_i$. 
See for example the four vectors 
\[
\begin{array}{c}
\textbf{\textit{x}}=(7,-7,\ 0,\ \ 0) \\
\textbf{\textit{y}}=(5,-5,\ 5,\ \ 0) \\
\textbf{w}=(3,-3,\ 3,-3)\\
\textbf{\textit{z}}=(2,-2,\ 2,-2).
\end{array}
\]
An aggregation function $G$ is said to be bipolar comonotone maxitive if it is maxitive on such a type of bipolar comonotonic \textit{bi-constants}, i.e. if fixed $K=\{1,\ldots,k\}$ it holds:
\begin{equation}
G\left(
{\bigvee_{i\in K}}^{b}{l_i\cdot \textbf{1}_{\left(A_i , B_i\right)}}
\right)
={\bigvee_{i\in K}}^b{G\left({l_i\cdot \textbf{1}_{\left(A_i , B_i\right)}}\right)}
.
\label{eq:bipmax}
\end{equation}
$G$ is said to be right bipolar comonotone maxitive if 
\begin{equation}
G\left(
{\bigvee_{i\in K}}^{b^+} {l_i\cdot \textbf{1}_{\left(A_i , B_i\right)}}
\right)
={\bigvee_{i\in K}}^{b^+}{G\left({l_i\cdot \textbf{1}_{\left(A_i , B_i\right)}}\right)}
.
\label{eq:posbipmax}
\end{equation}
$G$ is said to be left bipolar comonotone maxitive if 
\begin{equation}
G\left(
{\bigvee_{i\in K}}^{b^-} {l_i\cdot \textbf{1}_{\left(A_i , B_i\right)}}
\right)
={\bigvee_{i\in K}}^{b^-}{G\left({l_i\cdot \textbf{1}_{\left(A_i , B_i\right)}}\right)}
.
\label{eq:negbipmax}
\end{equation}
Clearly, due to bipolar comonotonicity, in equations \eqref{eq:bipmax}-\eqref{eq:negbipmax}: 
$$
{\bigvee_{i\in K}}^b {l_i\cdot \textbf{1}_{\left(A_i , B_i\right)}}=
{\bigvee_{i\in K}}^{b^+} {l_i\cdot \textbf{1}_{\left(A_i , B_i\right)}}
={\bigvee_{i\in K}}^{b^-} {l_i\cdot \textbf{1}_{\left(A_i , B_i\right)}}.
$$

\subsection{The bipolar Choquet integral}

\begin{Def} 
A function $\mu_b: \mathcal Q \rightarrow [-1,1]$ is a bi-capacity  \cite{grabisch2005biI,grabisch2005bi,greco2002bipolar}
on $N$  if
\begin{itemize}
	\item $\mu_b(\emptyset,  \emptyset)=0$,\  $\mu_b(N,  \emptyset)=1$ and  $\mu_b(\emptyset, N)=-1;$
		\item $\mu_b(A,B) \le \mu_b(C,D)\ \forall \  (A,B), (C,D) \in \mathcal Q$ such that  $(A,B)\precsim(C,D).$ 
\end{itemize}
\label{def:bi-capacity}
\end{Def}
\begin{Def}The bipolar Choquet integral of $\textbf{\textit{x}}=\left(x_1, \ldots, x_n\right)\in\left[-1,1\right]^n$ with respect to the bi-capacity $\mu_b$ is given by \cite{grabisch2005biI,grabisch2005bi,greco2002bipolar,greco2011choquet}:
\begin{equation}
Ch_b(\textbf{\textit{x}}, \mu_b) = \int_0^\infty \mu_b(\{i \in N:x_i > t \},\{i \in N:x_i  < -t \})dt.
\label{eq:bipcho}
\end{equation}
\label{Def:bipcho}
\end{Def}
\noindent The bipolar Choquet integral of $\textbf{\textit{x}}=\left(x_1, \ldots, x_n\right)\in\left[-1,1\right]^n$ with respect to the  bi-capacity $\mu_b$ can be rewritten as
\begin{equation}
Ch_b(\textbf{\textit{x}}, \mu_b) = \sum _{i=1}^n{\left(|x_{(i)}|-|x_{(i-1)}|\right)\mu_{b}(\{j\in N: x_j \ge |x_{(i)}|\},\{j\in N: x_j \le -|x_{(i)}|\})},
\label{eq:bipcho1}
\end{equation}
being $():N\rightarrow N$ any permutation of index such that $0=|x_{(0)}|\leq |x_{(1)}|\le\ldots\le |x_{(n)}|$. 
Note that to ensure that the pair
$\left(\{j\in N: x_j \ge |t|\}, \{j\in N: x_j \le -|t|\}\right)$ is an element of $\mathcal Q$ for all $t\in\rea$, we adopt the convention - which will be maintained trough all the paper - that in the case of $t=0$ the inequality $x_j \le -|0|=0$ must be intended as 
$x_j <0.$ 
The formulation \eqref{eq:bipcho1} will be useful in proving some results, like that exposed in the next representation theorem. 

\begin{Theo}
\cite{greco2002bipolar} An aggregation function $G: \left[-1,1\right]^n\rightarrow\left[-1,1\right]$ is idempotent and bipolar comonotonic additive if and only if there exists a bi-capacity $\mu_b$ such that, for all $\textbf{\textit{x}}\in \left[-1,1\right]^n$,
$$G(\textbf{\textit{x}})=Ch_b(\textbf{\textit{x}},\mu_b). $$
\label{bipcho}
\end{Theo}

\begin{rem}
Although the bipolar Choquet integral is trivially homogeneous, this condition does not appear in the theorem, since an aggregation function which is idempotent and bipolar comonotone additive is also homogeneous.  
Observe also that we could relax idempotency with the conditions $G(\textbf{1}_{(N,\emptyset)})=1$ and $G(\textbf{1}_{(\emptyset,N)})=-1$.
\end{rem}

\subsection{The bipolar Shilkret integral}

\begin{Def}The bipolar Shilkret integral of $\textbf{\textit{x}}=\left(x_1, \ldots, x_n\right)\in\left[-1,1\right]^n$ with respect to the bi-capacity $\mu_b$ is given by:
\begin{equation}
Sh_b(\textbf{\textit{x}},\mu_{b})={\bigvee_{i \in N}}^b\left\{|x_i|\cdot\mu_{b}(\{j\in N: x_j \ge |x_i|\},\{j\in N: x_j \le -|x_i|\})\right\}.
\label{eq:bipshi}
\end{equation}
\end{Def}

\begin{Def}The right bipolar Shilkret integral of $\textbf{\textit{x}}=\left(x_1, \ldots, x_n\right)\in\left[-1,1\right]^n$ with respect to the bi-capacity $\mu_b$ is given by:
\begin{equation}
Sh^+_b(\textbf{\textit{x}},\mu_{b})={\bigvee_{i \in N}}^{b^+}\left\{|x_i|\cdot\mu_{b}(\{j\in N: x_j \ge |x_i|\},\{j\in N: x_j \le -|x_i|\})\right\}.
\label{eq:posbipshi}
\end{equation}
\end{Def}

\begin{Def}The left bipolar Shilkret integral of $\textbf{\textit{x}}=\left(x_1, \ldots, x_n\right)\in\left[-1,1\right]^n$ with respect to the bi-capacity $\mu_b$ is given by:
\begin{equation}
Sh^-_b(\textbf{\textit{x}},\mu_{b})={\bigvee_{i \in N}}^{b^-}\left\{|x_i|\cdot\mu_{b}(\{j\in N: x_j \ge |x_i|\},\{j\in N: x_j \le -|x_i|\})\right\}.
\label{eq:negbipshi}
\end{equation}
\end{Def}

\noindent Clearly the three definitions are linked via the 
$$Sh_b(\textbf{\textit{x}},\mu_{b})={\bigvee}^b\left\{Sh^+_b(\textbf{\textit{x}},\mu_{b}), Sh^-_b(\textbf{\textit{x}},\mu_{b})\right\}.$$ 
The condition $Sh_b(\textbf{\textit{x}},\mu_{b})=0$ is equivalent to the $Sh^+_b(\textbf{\textit{x}},\mu_{b}) = - Sh^-_b(\textbf{\textit{x}},\mu_{b})$ and, in this case, or the three integrals  are all zero or they give three different results, one zero, one positive and one negative. 
We can think about them in terms of a neutral, an optimistic and a pessimistic aggregate evaluation of \textbf{\textit{x}}.
The condition $Sh_b(\textbf{\textit{x}},\mu_{b})\neq 0$ implies that $Sh^+_b(\textbf{\textit{x}},\mu_{b}) = Sh_b^-(\textbf{\textit{x}},\mu_{b})= Sh_b(\textbf{\textit{x}},\mu_{b})$.\\
The following theorems characterize the bipolar Shilkret integral. 

\begin{Theo}
An aggregation function $G: \left[-1,1\right]^n\rightarrow\left[-1,1\right]$ is idempotent, bipolar comonotone maxitive and homogeneous if and only if there exists a bi-capacity $\mu_{b}$ on $N$ such that, for all $\textbf{\textit{x}}\in \left[-1,1\right]^n$,
$$G(\textbf{\textit{x}})=Sh_b(\textbf{\textit{x}},\mu_{b}). $$
\label{bipshi}
\end{Theo}

\begin{rem}\label{rem11}
Let us note that theorem \ref{bipshi} implies, as corollary, theorem \ref{theo:shilkret} since bipolar comonotone maxitivity restricted on $\rea^+$ implies comonotone maxitivity.
\end{rem}

\begin{Theo}
An aggregation function $G: \left[-1,1\right]^n\rightarrow\left[-1,1\right]$ is idempotent, positive bipolar comonotone maxitive and homogeneous if and only if there exists a bi-capacity $\mu_{b}$ on $N$ such that, for all $\textbf{\textit{x}}\in \left[-1,1\right]^n$,
$$G(\textbf{\textit{x}})=Sh^+_b(\textbf{\textit{x}},\mu_{b}). $$
\label{posbipshi}
\end{Theo}

\begin{Theo}
An aggregation function $G: \left[-1,1\right]^n\rightarrow\left[-1,1\right]$ is idempotent, negative  bipolar comonotone maxitive and homogeneous if and only if there exists a bi-capacity $\mu_{b}$ on $N$ such that, for all $\textbf{\textit{x}}\in \left[-1,1\right]^n$,
$$G(\textbf{\textit{x}})=Sh^-_b(\textbf{\textit{x}},\mu_{b}). $$
\label{negbipshi}
\end{Theo}

\begin{rem}
Idempotency could be relaxed with the conditions $G(\textbf{1}_{(N,\emptyset)})=1$ and $G(\textbf{1}_{(\emptyset,N)})=-1$, in fact from these and from homogeneity idempotency can be elicited.
\end{rem}

\subsection{The bipolar Sugeno integral}

\begin{Def}
The bipolar Sugeno integral of a vector $\textbf{\textit{x}}=\left(x_1, \ldots, x_n\right)\in\left[-1,1\right]^n$ with respect to the bi-capacity $\mu_b$ on $N$ is given by:
\begin{eqnarray}
Su_b(\textbf{\textit{x}},\mu_b)&=&
{\bigvee_{i\in N}}^b\Bigl\{\textnormal{sign}\left(\mu_b\left(\{j\in N:x_j\ge |x_i|\},\{j\in N:x_j\le-|x_i|\}\right)\right)\cdot\nonumber\\
&&\qquad \cdot\bigwedge \left\{\left|\mu_b(\{j\in N:x_j\ge |x_i|\},\{j \in N:x_j\le-|x_i|\})\right|, |x_i|\right\}\Bigr\}.
\label{eq:bipsug}
\end{eqnarray}
\end{Def}

\begin{Def}
The right bipolar Sugeno integral of a vector $\textbf{\textit{x}}=\left(x_1, \ldots, x_n\right)\in\left[-1,1\right]^n$ with respect to the bi-capacity $\mu_b$ on $N$ is given by:
\begin{eqnarray}
Su^+_b(\textbf{\textit{x}},\mu_b)&=&
{\bigvee_{i\in N}}^{b^+}\Bigl\{\textnormal{sign}\left(\mu_b\left(\{j\in N:x_j\ge |x_i|\},\{j\in N:x_j\le-|x_i|\}\right)\right)\cdot\nonumber\\
&&\qquad \cdot\bigwedge \left\{\left|\mu_b(\{j\in N:x_j\ge |x_i|\},\{j \in N:x_j\le-|x_i|\})\right|, |x_i|\right\}\Bigr\}.
\label{eq:posbipsug}
\end{eqnarray}
\end{Def}

\begin{Def}
The left bipolar Sugeno integral of a vector $\textbf{\textit{x}}=\left(x_1, \ldots, x_n\right)\in\left[-1,1\right]^n$ with respect to the bi-capacity $\mu_b$ on $N$ is given by:
\begin{eqnarray}
Su^-_b(\textbf{\textit{x}},\mu_b)&=&
{\bigvee_{i\in N}}^{b^-}\Bigl\{\textnormal{sign}\left(\mu_b\left(\{j\in N:x_j\ge |x_i|\},\{j\in N:x_j\le-|x_i|\}\right)\right)\cdot\nonumber\\
&&\qquad \cdot\bigwedge \left\{\left|\mu_b(\{j\in N:x_j\ge |x_i|\},\{j \in N:x_j\le-|x_i|\})\right|, |x_i|\right\}\Bigr\}.
\label{eq:negbipsug}
\end{eqnarray}
\end{Def} 

\noindent Clearly the three definitions are linked via the 
$$Su_b(\textbf{\textit{x}},\mu_{b})={\bigvee}^b\left\{Su^+_b(\textbf{\textit{x}},\mu_{b}), Su^-_b(\textbf{\textit{x}},\mu_{b})\right\}.$$ 
The condition $Su_b(\textbf{\textit{x}},\mu_{b})=0$ is equivalent to the $Su_b^+(\textbf{\textit{x}},\mu_{b}) = - Su_b^-(\textbf{\textit{x}},\mu_{b})$ and, in this case, or the three integrals  are all zero or they give three different results, one zero (neutral), one positive (optimistic) and one negative (pessimistic). 
The condition $Su_b(\textbf{\textit{x}},\mu_{b})\neq 0$ implies that $Su_b^+(\textbf{\textit{x}},\mu_{b}) = Su_b^-(\textbf{\textit{x}},\mu_{b})= Su_b(\textbf{\textit{x}},\mu_{b})$.\\
The following theorems characterize the bipolar Sugeno integral. 

\begin{Theo}
An aggregation function $G: \left[-1,1\right]^n\rightarrow \left[-1,1\right]$ is idempotent, bipolar comonotone maxitive, bipolar stable with respect to the sign and bipolar stable with respect to the minimum if and only if there exists a  bi-capacity $\mu_b$ on $N$ such that, for all $\textbf{\textit{x}}\in \left[-1,1\right]^n$,
$$G(\textbf{\textit{x}})=Su_b(\textbf{\textit{x}},\mu_b). $$
\label{bipsug}
\end{Theo}

\begin{Theo}
An aggregation function $G: \left[-1,1\right]^n\rightarrow \left[-1,1\right]$ is idempotent, positive bipolar comonotone maxitive, bipolar stable with respect to the sign and bipolar stable with respect to the minimum if and only if there exists a  bi-capacity $\mu_b$ on $N$ such that, for all $\textbf{\textit{x}}\in \left[-1,1\right]^n$,
$$G(\textbf{\textit{x}})=Su_b^+(\textbf{\textit{x}},\mu_b). $$
\label{posbipsug}
\end{Theo}

\begin{Theo}
An aggregation function $G: \left[-1,1\right]^n\rightarrow \left[-1,1\right]$ is idempotent, negative bipolar comonotone maxitive, bipolar stable with respect to the sign and bipolar stable with respect to the minimum if and only if there exists a  bi-capacity $\mu_b$ on $N$ such that, for all $\textbf{\textit{x}}\in \left[-1,1\right]^n$,
$$G(\textbf{\textit{x}})=Su_b^-(\textbf{\textit{x}},\mu_b). $$
\label{negbipsug}
\end{Theo}

\section{Concluding remarks}\label{concluding}
In recent years there has been an increasing interest in development of new integrals useful in decision analysis process or in modeling engineering problems. 
An interesting line of research is that of bipolar fuzzy integrals, that considers the case in which the underling scale is bipolar. 
In this paper we have axiomatically characterized the bipolar Choquet integral and defined and axiomatically characterized the bipolar Shilkret integral and the bipolar Sugeno integral. 
Thus, the scenario of bipolar fuzzy integrals appears clearer and richer. 

\section{Appendix}

\textit{Proof of Theorem \ref{theo:shilkret}.}\\
First we prove the necessary part. 
Let us suppose there exists a capacity $\mu$ on $N$ such that, for all $\textbf{\textit{x}}\in \left(\alpha,\beta\right)^n$, $G(\textbf{\textit{x}})=Sh(\textbf{\textit{x}},\mu)$. 
In this case it is trivial to prove that the Shilkret integral is idempotent, comonotone maxitive and homogeneous by definition and we leave the proof to the reader. 
Now we prove the sufficient part of the theorem.
Let us define 
\begin{equation}
\mu(A)=G(\textbf{1}_A), \quad\textnormal{for all}\quad A\in 2^N.
\label{eq:capacity}
\end{equation}
Because $G$ is an idempotent aggregation function, we get $\mu(\emptyset)=0$, $\mu(N)=1$ and $\mu(A)\le \mu(B)$ whenever $A\subseteq B$. 
Thus $\mu$ is a capacity on $N$.
Every $\textbf{\textit{x}}=\left(x_1, \ldots, x_n\right)\in \left(\alpha, \beta\right)^n$ can be written as 
$$\textbf{\textit{x}}=
\bigvee_{i \in N}
x_{(i)}\cdot \textbf{1}_{\left\{j\in N\ |\ x_j\ge x_{(i)}\right\}}
$$
being $():N\rightarrow N$ any permutation of index such that $x_{(1)}\le\ldots\le x_{(n)}$.
Because vectors 
$
x_{(i)}\cdot \textbf{1}_{\left\{j\in N\ |\ x_j\ge x_{(i)}\right\}}
$
are comonotonic, we get the thesis by applying comonotonic maxitivity, homogeneity of $G$ and the definition of $\mu$ according to \eqref{eq:capacity}:
\[
G(\textbf{\textit{x}})=
G\left(
\bigvee_{i \in N}
x_{(i)}\cdot \textbf{1}_{\left\{j\in N\ |\ x_j\ge x_{(i)}\right\}}
\right)=
\bigvee_{i \in N}
G\left(
x_{(i)}\cdot \textbf{1}_{\left\{j\in N\ |\ x_j\ge x_{(i)}\right\}}
\right)=
\]
\[
=\bigvee_{i \in N}
x_{(i)}\cdot 
G\left(
\textbf{1}_{\left\{j\in N\ |\ x_j\ge x_{(i)}\right\}}
\right)=
\bigvee_{i \in N}
x_{(i)}\cdot \mu\left(
\left\{j\in N\ |\ x_j\ge x_{(i)}\right\}
\right)=
Sh(\textbf{\textit{x}},\mu)
\]
\begin{flushright}
$\square$
\end{flushright}

\noindent\textit{Proof of Theorem \ref{bipcho}.}\\
First we prove the necessary part.
Let us suppose that there exists a bi-capacity $\mu_b$ such that, for all $\textbf{\textit{x}}\in \left[-1,1\right]^n$,
$G(\textbf{\textit{x}})=Ch_b(\textbf{\textit{x}},\mu_b).$ 
Idempotency of the bipolar Choquet integral follows from definition, because if 
$\lambda\geq 0$, then 
$Ch_b\left(\lambda\cdot\textbf{1}_{\left(N,\emptyset\right)},\mu_b\right)
=\int_0^\lambda{\mu_b\left(N,\emptyset\right)dt}=\lambda$, while if  
$\lambda< 0$, then 
$Ch_b\left(\lambda\cdot\textbf{1}_{\left(N,\emptyset\right)},\mu_b\right)=
\int_0^{-\lambda}{\mu_b\left(\emptyset,N\right)dt}=\lambda$.
If $\textbf{\textit{x}}$, $\textbf{\textit{y}}\in \left[-1,1\right]^n$ are bipolar comonotone, 
then there exists a permutation of indexes $():N\rightarrow N$ such that 
$0=|x_{(0)}|\leq |x_{(1)}|\le\ldots\le |x_{(n)}|$ and $0=|y_{(0)}|\leq |y_{(1)}|\le\ldots\le |y_{(n)}|$, and thus   
\[
Ch_b\left(\textbf{\textit{x}},\mu_b\right)= \sum _{i=1}^n{\left(|x_{(i)}|-|x_{(i-1)}|\right)\cdot \mu_b\left(\{j\in N: x_j \ge |x_{(i)}|\},\{j\in N: x_j \le -|x_{(i)}|\}\right)},
\]
and 
\[
Ch_b\left(\textbf{\textit{y}},\mu_b\right)= \sum _{i=1}^n{\left(|y_{(i)}|-|y_{(i-1)}|\right)\cdot \mu_b\left(\{j\in N: y_j \ge |y_{(i)}|\},\{j\in N: y_j \le -|y_{(i)}|\}\right)}.
\]
Since $\textbf{\textit{x}}$ and $\textbf{\textit{y}}$ are absolutely comonotonic an cosigned, for every $i=1,\ldots,n$
\begin{equation}
\mu_b\left(\{j\in N: x_j \ge |x_{(i)}|\},\{j\in N: x_j \le -|x_{(i)}|\}\right)=
\mu_b\left(\{j\in N: y_j \ge |y_{(i)}|\},\{j\in N: y_j \le -|y_{(i)}|\}\right).
\label{eq:n}
\end{equation}
Moreover, again because $\textbf{\textit{x}}$ and $\textbf{\textit{y}}$ are absolutely comonotonic and cosigned, for every $i=1,\ldots,n$, 
$|x_{(i)} + y_{(i)}|= |x_{(i)}| + |y_{(i)}|$ and consequently 
\begin{equation}
0=|x_{(0)}+y_{(0)}|\leq |x_{(1)}+ y_{(i)}|\le\ldots\le |x_{(n)}+y_{(n)}|
\quad
\textnormal{for every}
\quad 
i=1,\ldots,n.
\label{eq:nn}
\end{equation}
By \eqref{eq:n} and \eqref{eq:nn} we get 
$Ch_b\left(\textbf{\textit{x}},\mu_b\right)
+Ch_b\left(\textbf{\textit{y}},\mu_b\right)=
Ch_b\left(\textbf{\textit{x}}+\textbf{\textit{y}},\mu_b\right)$.\\
Now we prove the sufficient part of the theorem. 
Let us define  
\begin{equation}
\mu_b(A,B)=G\left(\textbf{1}_{(A,B)}\right), \quad\textnormal{for all}\quad (A,B)\in\mathcal Q.
\label{eq:bicapacity}
\end{equation} 
$\mu_b$ represents a bi-capacity, because by idempotency of $G$ we get that 
$
\mu_b(N,\emptyset)=G\left(\textbf{1}_{(N,\emptyset)}\right)=1, 
$ 
$
\mu_b(\emptyset,N)=G\left(\textbf{1}_{(\emptyset,N)}\right)=-1,
$ 
$
\mu_b(\emptyset,\emptyset)=G\left(\textbf{1}_{(\emptyset,\emptyset)}\right)=0
$. 
Moreover, if 
$(A,B)\precsim (A',B')$,  
being for all $i\in N$, the $i^{th}$ component of the vector 
$
\textbf{1}_{(A,B)}
$
not greater than the $i^{th}$ component of the vector 
$\textbf{1}_{(A',B')}$
and being $G$ an aggregation function (then monotone), thus  
$\mu_b(A,B)\le\mu_b(A',B').$
Observe now that any vector $\textbf{\textit{x}}=\left(x_1, \ldots, x_n\right)\in\left[-1,1\right]^n$ can be rewritten as 
\begin{equation}
\textbf{\textit{x}}= \sum _{i=1}^n{\left(|x_{(i)}|-|x_{(i-1)}|\right)\cdot \textbf{1}_{\left(\{j\in N: x_j \ge |x_{(i)}|\},\{j\in N: x_j \le -|x_{(i)}|\}\right)}},
\label{eq:dec}
\end{equation}
being $():N\rightarrow N$ any permutation of indexes such that $0=|x_{(0)}|\leq |x_{(1)}|\le\ldots\le |x_{(n)}|.$ 
Let us note that for all $(A,B),\ (A',B')\in \mathcal Q$ such that $(A,B)\subseteq (A',B')$ and for all  $a,b\in [0,1]$, vectors $a\cdot \textbf{1}_{(A,B)}$ and $b\cdot \textbf{1}_{(A',B')}$ are bipolar comonotone. 
Consequently, \eqref{eq:dec} shows that any vector $\textbf{\textit{x}}\in\left[-1,1\right]^n$ can be decomposed as a sum of bipolar comonotonic vectors. 
Remembering that an aggregation function which is idempotent and bipolar comonotone additive is also homogeneous, 
thus to get the thesis it is sufficient to apply, respectively, bipolar comonotone additivity, homogeneity of $G$ and definition of bi-capacity $\mu_b$ according to \eqref{eq:bicapacity}: 
\[
G(\textbf{\textit{x}}) =G\left(
\sum _{i=1}^n{\left(|x_{(i)}|-|x_{(i-1)}|\right)\cdot \textbf{1}_{\left(\{j\in N: x_j \ge |x_{(i)}|\},\{j\in N: x_j \le -|x_{(i)}|\}\right)}}
\right)
=
\]
\[
=\sum _{i=1}^n{\left(|x_{(i)}|-|x_{(i-1)}|\right)\cdot G\left(\textbf{1}_{\left(\{j\in N: x_j \ge |x_{(i)}|\},\{j\in N: x_j \le -|x_{(i)}|\}\right)}\right)}
=Ch_b(\textbf{\textit{x}}, \mu_b).
\]
\begin{flushright}
$\square$
\end{flushright}

\noindent \textit{Proof of Theorem \ref{bipshi}.}\\
First we prove the necessary part.
Let us suppose there exists a bi-capacity $\mu_b$ such that, for all $\textbf{\textit{x}}\in \left[-1,1\right]^n$,
$G(\textbf{\textit{x}})=Sh_b(\textbf{\textit{x}},\mu_b).$ 
The bipolar Shilkret integral is, trivially, idempotent and homogeneous and we only need to demonstrate the bipolar comonotonic maxitivity.
Let us consider a set of indexes $K=\{1,\ldots,k\}$, $k$ increasing levels $l_1,\ldots,l_k\in\rea$ with
$0<l_1<l_2<\ldots<l_k\le 1$ 
and a sequence 
$\left\{\left(A_i , B_i\right)\right\}_{i\in K}$ 
such that
$\left(A_i , B_i\right)\in\mathcal Q$ and  
$\left(A_{i+1} , B_{i+1}\right)\subseteq \left(A_i , B_i\right)$ 
for all $i\in K.$
The $j^{th}$ component of the vector 
$
\bigvee_{i\in K}^b\{l_i\cdot \textbf{1}_{\left(A_i , B_i\right)}\}
$
is equal to $l_i$ if $j\in A_i\setminus A_{i+1}$, 
is equal to $-l_i$ if $j\in B_i\setminus B_{i+1}$ and is equal to zero if $j\in N\setminus \left(A_1\cup B_1\right)$ for all $i\in K$ and taking $A_{k+1}=B_{k+1}=\emptyset$. 
Clearly, such a vector has a component greater or equal to $l_i$ for indexes in  
$A_i$ 
and has component smaller or equal to $-l_i$ for indexes in  
$B_i$. 
Thus, by definition
\begin{equation}
Sh_b\left(
{\bigvee_{i\in K}}^{b}\{l_i\cdot \textbf{1}_{\left(A_i , B_i\right)}\}
, \mu_b
\right)
=
{\bigvee_{i\in K}}^b
\left\{
l_i\cdot \mu_b\left(\left(A_i,B_i\right)\right)
\right\}
={\bigvee_{i\in K}}^b\{Sh_b\left({l_i\cdot \textbf{1}_{\left(A_i , B_i\right)}},\mu_b\right)\}
.
\label{eq:bipmaxproof}
\end{equation}
Now we prove the sufficient part of the theorem. 
Let us define 
\begin{equation}
\mu_b(A,B)=G\left(\textbf{1}_{(A,B)}\right), \quad\textnormal{for all}\quad (A,B)\in\mathcal Q.
\label{eq:bicapacity1}
\end{equation} 
$\mu_b$ represents a bi-capacity (see proof of theorem \ref{bipcho}). 
Notice that each $\textbf{\textit{x}}\in\left[-1,1\right]^n$ can be rewritten as
\begin{equation}
\textbf{\textit{x}}={\bigvee_{i\in N}}^b|x_{i}|
\cdot
 \textbf{1}_{
 \left(\left\{j\ |\ x_j\ge |x_{i}|\right\},\left\{j\ |\ x_j\le -|x_{i}|\right\}\right)}
\label{eq:11}
\end{equation}
and observe that vectors  
$|x_{i}|\cdot \textbf{1}_{
 \left(\left\{j\in N\  |\ x_j\ge |x_{i}|\right\},\left\{j\in N\ |\ x_j\le -|x_{i}|\right\}\right)}$, $i=1\ldots,n$  
are bipolar comonotone. 
Consequently, for any $\textbf{\textit{x}}\in\left[-1,1\right]^n$ by bipolar comonotone maxitivity, homogeneity and definition of bi-capacity $\mu_b$ according to the \eqref{eq:bicapacity1} we get
\[
G(\textbf{\textit{x}})=G\left({\bigvee_{i\in N}}^b|x_{i}|
\cdot
 \textbf{1}_{
 \left(\left\{j\ |\ x_j\ge |x_{i}|\right\},\left\{j\ |\ x_j\le -|x_{i}|\right\}\right)}\right)
={\bigvee_{i\in N}}^bG\left(|x_{i}|
\cdot
 \textbf{1}_{
 \left(\left\{j\ |\ x_j\ge |x_{i}|\right\},\left\{j\ |\ x_j\le -|x_{i}|\right\}\right)}\right)=
\]
\[
={\bigvee_{i\in N}}^b|x_{i}|
\cdot G\left(
 \textbf{1}_{
 \left(\left\{j\ |\ x_j\ge |x_{i}|\right\},\left\{j\ |\ x_j\le -|x_{i}|\right\}\right)}\right)
 =
 {\bigvee_{i\in N}}^b|x_{i}|
\cdot 
\mu_b
\left(
\left\{j\ |\ x_j\ge |x_{i}|\right\},\left\{j\ |\ x_j\le -|x_{i}|\right\}
\right)=Sh_b(\textbf{\textit{x}},\mu_b)
\]
\begin{flushright}
$\square$
\end{flushright}

\noindent \textit{Proof of Theorems \ref{posbipshi} and \ref{negbipshi}.} 
They are analogous to the proof of previous \textit{Theorem} \ref{bipshi}.
\begin{flushright}
$\square$
\end{flushright}

\noindent \textit{Proof of Theorem \ref{bipsug}.}
First we prove the necessary part.
Let us suppose there exists a bi-capacity $\mu_b$ such that, for all $\textbf{\textit{x}}\in \left[-1,1\right]^n$,
$G(\textbf{\textit{x}})=Su_b(\textbf{\textit{x}},\mu_b).$ 
The Sugeno integral is idempotent by definition.
Bipolar stability with respect to the sign and with respect to the minimum are trivially verified once we consider that for all $r>0$ and for all $(A,B)\in \mathcal Q$
\[
Su_b
\left(
r\cdot
\textbf{1}_{(A,B)},\mu_b
\right)=
sign\left(\mu_b(A,B)\right)
\bigwedge
\left\{
r,
\left|
\mu_b(A,B)
\right|
\right\}
.
\]
Let us consider a set of indexes $K=\{1,\ldots,k\}$, $k$ increasing levels $l_1,\ldots,l_k\in\rea$ with
$0<l_1<l_2<\ldots<l_k\le 1$
and a sequence 
$\left\{\left(A_i , B_i\right)\right\}_{i\in K}$ 
such that
$\left(A_i , B_i\right)\in\mathcal Q$ and  
$\left(A_{i+1} , B_{i+1}\right)\subseteq \left(A_i , B_i\right)$ 
for all $i\in K.$
Thus, by definition
\begin{eqnarray}
Su_b\left(
{\bigvee_{i\in K}}^{b}\{l_i\cdot \textbf{1}_{\left(A_i , B_i\right)}\}
, \mu_b
\right)&
=&
{\bigvee_{i\in K}}^b
\left\{
sign 
\left[
\mu_b\left(\left(A_i,B_i\right)\right)
\right]
\bigwedge
\left\{
l_i, 
|\mu_b\left(\left(A_i,B_i\right)\right)|
\right\}
\right\}=\nonumber\\
&=&{\bigvee_{i\in K}}^b\{Su_b\left({l_i\cdot \textbf{1}_{\left(A_i , B_i\right)}},\mu_b\right)\}
.
\label{eq:bipmaxproof1}
\end{eqnarray}
Now we prove the sufficient part of the theorem. 
Let us define  
$
\mu_b(A,B)=G\left(\textbf{1}_{(A,B)}\right)
$
for all $(A,B)\in\mathcal Q$. 
$\mu_b$ represents a bi-capacity (see proof of theorem \ref{bipcho}). 
Let us note that using bipolar stability with respect to the minimum and idempotency of $G$ we have that for all $r>0$ and for all $(A,B)\in\mathcal Q$, 
\begin{equation}
\left|
G\left(r\cdot \textbf{1}_{(A,B)}\right)\right|
=\bigwedge\left\{r,\ \left|G\left(\textbf{1}_{(A,B)}\right)\right|\right\}.
\label{eq:po}
\end{equation}
The \eqref{eq:po} is obvious if $r=0$ or $r=1$. 
If $0<r<1$ and 
$\left|G\left(\textbf{1}_{(A,B)}\right)\right|>\left|G\left(r\cdot \textbf{1}_{(A,B)}\right)\right|$, then 
using stability with respect to the minimum, 
$\left|G\left(r\cdot \textbf{1}_{(A,B)}\right)\right|=r $ and the \eqref{eq:po} is true again. 
If 
$\left|G\left(\textbf{1}_{(A,B)}\right)\right|=\left|G\left(r\cdot \textbf{1}_{(A,B)}\right)\right|$ 
observe that by monotonicity and idempotency of $G$, 
$\left|G\left(r\cdot \textbf{1}_{(A,B)}\right)\right|\leq \left|G\left(r\cdot \textbf{1}_{(N,\emptyset)}\right)\right|=r$,
which means that also in this last case the \eqref{eq:po} is true.
Finally, notice that each $\textbf{\textit{x}}\in\left[-1,1\right]^n$ can be rewritten as
\begin{equation}
\textbf{\textit{x}}={\bigvee_{i\in N}}^b|x_{i}|
\cdot
 \textbf{1}_{
 \left(\left\{j\ |\ x_j\ge |x_{i}|\right\},\left\{j\ |\ x_j\le -|x_{i}|\right\}\right)}
\label{eq:111}
\end{equation}
and observe that vectors  
$|x_{i}|\cdot \textbf{1}_{
 \left(\left\{j\in N\  |\ x_j\ge |x_{i}|\right\},\left\{j\in N\ |\ x_j\le -|x_{i}|\right\}\right)}$, $i=1\ldots,n$  
are bipolar comonotone.\\
Consequently, for any $\textbf{\textit{x}}\in\left[-1,1\right]^n$ by bipolar comonotone maxitivity 
\[
G(\textbf{\textit{x}})=G\left({\bigvee_{i\in N}}^b|x_{i}|
\cdot
 \textbf{1}_{
 \left(\left\{j\ |\ x_j\ge |x_{i}|\right\},\left\{j\ |\ x_j\le -|x_{i}|\right\}\right)}\right)
={\bigvee_{i\in N}}^b 
G\left(|x_{i}|
\cdot
 \textbf{1}_{
 \left(\left\{j\ |\ x_j\ge |x_{i}|\right\},\left\{j\ |\ x_j\le -|x_{i}|\right\}\right)}\right)=
\]
( by bipolar stability with respect to the sign )
\[
={\bigvee_{i\in N}}^b 
\left\{
\textnormal{sign}\left[
G\left(\textbf{1}_{ \left(\left\{j\ |\ x_j\ge |x_{i}|\right\},\left\{j\ |\ x_j\le-|x_{i}|\right\}\right)}\right)\right]
\left|G\left(|x_{i}|
\cdot
 \textbf{1}_{
 \left(\left\{j\ |\ x_j\ge |x_{i}|\right\},\left\{j\ |\ x_j\le -|x_{i}|\right\}\right)}\right)\right|\right\}=
\]
\[
={\bigvee_{i\in N}}^b 
\left\{
\textnormal{sign}
\left[
\mu_b\left(
\left\{j\ |\ x_j\ge |x_{i}|\right\},\left\{j\ |\ x_j\le-|x_{i}|\right\}
\right)
\right]
\left|G\left(|x_{i}|
\cdot
 \textbf{1}_{
 \left(\left\{j\ |\ x_j\ge |x_{i}|\right\},\left\{j\ |\ x_j\le -|x_{i}|\right\}\right)}\right)\right|
\right\}=
\]
( by bipolar stability with respect to the minimum )
\[
={\bigvee_{i\in N}}^b
\left\{
\textnormal{sign}
\left[
\mu_b\left(
\left\{j\ |\ x_j\ge |x_{i}|\right\},\left\{j\ |\ x_j\le-|x_{i}|\right\}
\right)
\right]
\bigwedge\left\{|x_{i}|,\ 
\left|
G\left(
 \textbf{1}_{
 \left(\left\{j\ |\ x_j\ge |x_{i}|\right\},\left\{j\ |\ x_j\le -|x_{i}|\right\}\right)}
\right)
\right|
\right\}
\right\}=
\]
\[
={\bigvee_{i\in N}}^b
\left\{
\textnormal{sign}
\left[
\mu_b\left(
\left\{j\ |\ x_j\ge |x_{i}|\right\},\left\{j\ |\ x_j\le-|x_{i}|\right\}
\right)
\right]
\bigwedge\left\{|x_{i}|,\ 
\left|
\mu_b\left(
\left\{j\ |\ x_j\ge |x_{i}|\right\},\left\{j\ |\ x_j\le-|x_{i}|\right\}
\right)
\right|
\right\}
\right\}
\]
that is the Sugeno integral $Su_b(\textbf{\textit{x}},\mu_b)$. 
\begin{flushright}
$\square$
\end{flushright}

\noindent \textit{Proof of Theorems \ref{posbipsug} and \ref{negbipsug}.} 
They are analogous to the proof of previous \textit{Theorem} \ref{bipsug}.
\begin{flushright}
$\square$
\end{flushright}



\begin{thebibliography}{10}

\bibitem{bodjanova2009sugeno}
S.~Bodjanova and M.~Kalina.
\newblock {Sugeno and Shilkret integrals, and T-and S-evaluators}.
\newblock In {\em SISY 2009 - 7th International Symposium on Intelligent
  Systems and Informatics}, pages 109--114.

\bibitem{choquet1953theory}
G.~Choquet.
\newblock {Theory of Capacities}.
\newblock {\em Annales de l'Institute Fourier (Grenoble)}, 5:131--295, 1953/54.

\bibitem{figueira2005multiple}
J.~Figueira, S.~Greco, and M.~Ehrgott.
\newblock {\em {Multiple criteria decision analysis: state of the art
  surveys}}, volume~78.
\newblock Springer Verlag, 2005.

\bibitem{grabisch1996application}
M.~Grabisch.
\newblock {The application of fuzzy integrals in multicriteria decision
  making}.
\newblock {\em European Journal of Operational Research}, 89(3):445--456, 1996.

\bibitem{grabisch2003symmetric}
M.~Grabisch.
\newblock {The symmetric Sugeno integral}.
\newblock {\em Fuzzy Sets and Systems}, 139(3):473--490, 2003.

\bibitem{grabisch2004mobius}
M.~Grabisch.
\newblock {The M{\"o}bius transform on symmetric ordered structures and its
  application to capacities on finite sets}.
\newblock {\em Discrete Mathematics}, 287(1-3):17--34, 2004.

\bibitem{grabisch2005biI}
M.~Grabisch and C.~Labreuche.
\newblock {Bi-capacities--I: definition, M{\"o}bius transform and interaction}.
\newblock {\em Fuzzy Sets and Systems}, 151(2):211--236, 2005.

\bibitem{grabisch2005bi}
M.~Grabisch and C.~Labreuche.
\newblock {Bi-capacities--II: the Choquet integral}.
\newblock {\em Fuzzy Sets and Systems}, 151(2):237--259, 2005.

\bibitem{grabisch2005fuzzy}
M.~Grabisch and C.~Labreuche.
\newblock Fuzzy measures and integrals in {MCDA}.
\newblock {\em Multiple criteria decision analysis: state of the art surveys},
  pages 563--604, 2005.

\bibitem{grabisch2009aggregation}
M.~Grabisch, J.L. Marichal, R.~Mesiar, and E.~Pap.
\newblock {\em {Aggregation Functions (Encyclopedia of Mathematics and its
  Applications)}}.
\newblock Cambridge University Press, 2009.

\bibitem{greco2011linz}
S.~Greco.
\newblock {Generalizing again the Choquet integral: the profile dependent
  Choquet integral}.
\newblock In Mesiar R. Klement E.~P. Dubois~D., Grabisch~M., editor, {\em
  Decision Theory: Qualitative and Quantitative Approaches}, pages 66--79. Linz
  Seminar on Fuzzy Set and System, 2011.

\bibitem{greco2011choquet}
S.~Greco, B.~Matarazzo, and S.~Giove.
\newblock {The Choquet integral with respect to a level dependent capacity}.
\newblock {\em Fuzzy Sets and Systems}, 175(1), 2011.

\bibitem{greco2002bipolar}
S.~Greco, B.~Matarazzo, and R.~Slowinski.
\newblock {Bipolar Sugeno and Choquet integrals}.
\newblock In G.~Pasi B.~De~Baets, J.~Fodor, editor, {\em EUROWorking Group on
  Fuzzy Sets}, pages 191--196. Workshop on Informations Systems (EUROFUSE
  2002), Varenna, Italy, 2002.

\bibitem{klement2010universal}
E.P. Klement, R.~Mesiar, and E.~Pap.
\newblock {A universal integral as common frame for Choquet and Sugeno
  integral}.
\newblock {\em IEEE Transactions on Fuzzy Systems}, 18(1):178--187, 2010.

\bibitem{marichal2001axiomatic}
J.L. Marichal.
\newblock {An axiomatic approach of the discrete Sugeno integral as a tool to
  aggregate interacting criteria in a qualitative framework}.
\newblock {\em IEEE Transactions on Fuzzy Systems}, 9(1):164--172, 2001.

\bibitem{mesiar2009level}
R.~Mesiar, A.~Mesiarov{\'a}-Zem{\'a}nkov{\'a}, and K.~Ahmad.
\newblock {Level-dependent Sugeno integral}.
\newblock {\em IEEE Transactions on Fuzzy Systems}, 17(1):167--172, 2009.

\bibitem{mesiar2010discrete}
R.~Mesiar, A.~Mesiarov{\'a}-Zem{\'a}nkov{\'a}, and K.~Ahmad.
\newblock {Discrete Choquet integral and some of its symmetric extensions}.
\newblock {\em Fuzzy Sets and Systems}, 184(1):148--155, 2011.

\bibitem{schmeidler1986integral}
D.~Schmeidler.
\newblock {Integral representation without additivity}.
\newblock {\em Proceedings of the American Mathematical Society},
  97(2):255--261, 1986.

\bibitem{shilkret1971maxitive}
N.~Shilkret.
\newblock {Maxitive measure and integration}.
\newblock In {\em Indagationes Mathematicae (Proceedings)}, volume~74, pages
  109--116. Elsevier, 1971.

\bibitem{sipos1979integral}
J.~{\v{S}}ipo{\v{s}}.
\newblock Integral with respect to a pre-measure.
\newblock {\em Mathematica Slovaca}, 29(2):141--155, 1979.

\bibitem{sugeno1974theory}
M.~Sugeno.
\newblock {\em {Theory of fuzzy integrals and its applications}}.
\newblock Ph.D. Thesis, Tokyo institute of Technology, 1974.
\end{thebibliography}
\end{document}